# SUBSTITUTION MARKOV CHAINS AND MARTIN BOUNDARIES

DAVID KOSLICKI[1*], MANFRED DENKER[2]

ABSTRACT. Substitution Markov chains have been introduced in [7] as a new model to describe molecular evolution. In this note, we study the associated Martin boundaries from a probabilistic and topological viewpoint. An example is given that, although having a boundary homeomorphic to the well-known coin tossing process, has a metric description that differs significantly.

**1. Introduction.** Deterministic substitutions are used to construct symbolic dynamical systems by means of iterating a string re-write rule on the letters of a finite alphabet. A classic example is the Thue-Morse sequence [12]

$$0110100110010110\cdots$$

which can be obtained by iterating the substitution $0 \to 01$, $1 \to 10$ repeatedly on the letter 0. In general, a substitution is defined as follows: We let $\mathcal{A}$ be a finite alphabet and let $\mathcal{A}^*$ denote the set of all finite length words formed via concatenation from $\mathcal{A}$. Then a *substitution* $\sigma$ is a morphism from $\mathcal{A}$ into the set $\mathcal{A}^*$ whose domain is then extended to $\mathcal{A}^*$ by concatenation: $\sigma(WV) = \sigma(W)\sigma(V)$. The symbolic dynamical system that arises from this substitution is then the orbit closure under the shift map on a fixed point obtained from iterating the substitution.

We wish to study randomized versions of such substitutions. A first analysis of such a randomization was performed by Peyriére ([8], [9], [10], [11]) in a series of articles where a construction of Mandelbrot (referred to by Peyriére as an *M-system*, or a *random beadset*, or a *random substitution*) was taken as the motivation. These articles give a definition of a random substitution and investigate the convergence of the frequencies of subwords under various conditions (they also investigate an associated kind of graph different that considered here). The main approach used in these papers is to use Doob's martingale convergence theorem or an estimation of second moments to obtain the desired results.

After Peyriére's papers, the concept of random substitutions lay dormant until Burton and Wing investigated the behavior of the complexity function for random substitutions [14].

We take as our starting point a generalization of Peyriére's definition of a random substitution, first introduced in the PhD thesis of one of the authors [7]. The main purpose of this note is to show how Martin boundaries can be determined for transient random substitution Markov chains in a probabilistic, topological and metric fashion. We do this by example, in which the key property is to replace the root $a$ with words $aw$ or $wa$ randomly. This motivates a large class of similar SMC's which are persistent and expanding (see definition 2.2) and where similar analysis as in this note can be performed. For example, by

2000 AMS *Mathematics subject classification.* 60J50, 31C35, 28D05, 37A50.
The authors were supported by the National Science Foundation under grant No. DMS-1008538.
*Corresponding Author: `david.koslicki@math.oregonstate.edu`.
Received by the editors March 17, 2015.

replacing $a$ with $a$, $aw$, or $wa$ (see the examples in section 3). For the example considered herein, the Martin boundary can be calculated explicitly and we analyze its metric properties in more detail (see Theorems 4.5 and 4.6).

**2. Random Substitutions and their Properties.** Let $\mathcal{A} = \{a_1, a_2, \ldots, a_N\}$ be an alphabet of $N$ symbols. For $\ell \geq 1$, let $\mathcal{A}^\ell$ be the set of all $\ell$-length words on $\mathcal{A}$, that is, all concatenations of $\ell$ symbols from $\mathcal{A}$. For $w \in \mathcal{A}^\ell$, we denote the length of $w$ as $|w| = \ell$. Denote the set of all finite length sequences as $\mathcal{A}^* = \cup_{\ell > 0} \mathcal{A}^\ell$. By a subword of $w$, we mean a word $u$ consisting of a contiguous sequence of symbols from $w$ and denote this relationship by $u \subseteq w$. Let $|w|_u$ denote the number of times $u$ appears as a subword of $w$. That is, $|w|_u$ is the cardinality of the set $w(u) = \{i \ : \ w_i w_{i+1} \ldots w_{i+|u|-1} = u\}$.

**Definition 2.1** (Substitution Markov Chain). *A substitution Markov chain is a Markov chain with state space $\mathcal{A}^*$ and transition operator $P$ satisfying the following property: For each $a \in \mathcal{A}$, there is a probability $P_a : \mathcal{A}^* \to [0, 1]$ such that for $v \in \mathcal{A}^*$, and $w = w_1 \ldots w_\ell$ with $w_i \in \mathcal{A}$,*

$$(2.1) \qquad P(w, v) = \sum_{\substack{v_1 \ldots v_\ell = v \\ v_i \in \mathcal{A}^*}} \prod_{i=1}^{\ell} P_{w_i}(v_i).$$

We will abbreviate "substitution Markov chain" by SMC. According to the above definition, to define an SMC, we need only specify the functions $P_a$ for each $a \in \mathcal{A}$. This definition appears in [7] for the first time. It generalizes Peyriére's definition and the notion of an S-adic transformation. Example 1 in section 3 is an SMC where the transition probabilities cannot be defined by Peyriére's definition of a random substitution (this follows by comparing our definition with [10]).

Substitution Markov chains form a subclass of all possible Markov chains with countable state space, and this class certainly contains all finite state Markov chains. It is worthwhile mentioning that an SMC defines a multitype branching process in a canonical way (see [2] for a definition). Let $\mathcal{A}$ denote the different types of a branching process. Then starting with one species $i$ of type $b \in \mathcal{A}$, the SMC creates words and we consider the letters appearing in it as descendants of $i$. If $X_i = (X_{i,a})_{a \in \mathcal{A}}$ denotes the vector of the number of descendants of types $a \in \mathcal{A}$, we can define a multitype branching process by

$$Z_0 = (1, 0, \ldots, 0) \qquad Z_n = \left( \sum_{i=1}^{Z_{n-1,a}} X_{i,a} \right)_{a \in \mathcal{A}},$$

where $Z_{n-1,a}$ denotes the $a$-th coordinate of $Z_{n-1}$. Clearly, $Z_n$ is a Markov chain as well, and its properties, in particular its Martin boundary, is well examined in the literature (see [1] and [4] for example).

It is clear from this that the investigation of SMCs and their boundary properties is more elaborate than similar questions for multitype branching processes. There is another difference to be mentioned here. In view of the application, we are not interested in studying the number of different genotypes/subletters, but rather the structure of genomes/words. This requires to consider topological and metric properties as well, which is usually not done for branching processes. In particular we are interested in the metrically isomorphic embedding of the state space into a compact space keeping the probabilistic structure meaningful. This viewpoint is in accordance with Martin's original idea.

We therefore define

**Definition 2.2.** *Let $(\mathcal{A}^*, (P_a)_{a \in \mathcal{A}})$ be a SMC.*

(1) *A letter $a \in \mathcal{A}$ is called **persistent** if there exists $w \in \mathcal{A}^*$ with $P_a(w) > 0$ and $|w|_a > 0$, and it is called **expanding** if in addition $|w| \geq 2$.*
(2) *A substitution Markov chain is called **persistent** if each $a \in \mathcal{A}$ is persistent.*
(3) *An element $a \in \mathcal{A}$ is called a **root** of the Markov chain if for each $w \in \mathcal{A}^*$ with $P_a(w) > 0$ it follows that $|w|_a = 1$ and for each $a \neq b \in \mathcal{A}^*$ and $w \in \mathcal{A}^*$ with $P_b(w) > 0$, $a$ is not a subword of $w$: $a \not\subset w$.*
(4) *An SMC is called **constant length** if there exists an $L \in \mathbb{N}$ such that for all $a \in \mathcal{A}$ and $w \in \mathcal{A}^*$ with $P_a(w) > 0$, $|w| = L$.*

**Lemma 2.3.** *Let $a \in \mathcal{A}$ be a root for the SMC $(\mathcal{A}^*, (P_c)_{c \in \mathcal{A}})$. Then every word $w \in \mathcal{A}^*$ with $n$-step transition probabilities $P^{(n)}(a, w) > 0$ for some $n \geq 0$ has a unique representation as a subset in $\mathcal{A}^{\mathbb{Z}}$:*
$$[w] = \{(x_k)_{k \in \mathbb{Z}} : x_k = w_k \in \mathcal{A}, r \leq k \leq s\}$$
*where $w_0 = a$ and $w_r w_{r+1} ... w_s = w$.*

*Proof.* By induction one shows that each $w$ with $P^{(n)}(a, w) > 0$ contains the letter $a$ exactly once. □

**Lemma 2.4.** *A substitution Markov chain $(\mathcal{A}^*, (P_c)_{c \in \mathcal{A}})$ with a persistent and expanding letter $a \in \mathcal{A}$ is transient under the initial distribution $\epsilon_a(\cdot)$, the point mass in $a$.*

*Proof.* Let $\eta = \sum_{w \in \mathcal{A}^*: |w| \geq 2} P_a(w)$. Since $P(v, w) > 0$ implies that $|w| \geq |v|$ we must have by persistence that
$$P^{(n)}(v, v) \leq (1 - \eta)^n$$
hence $v$ returns to itself only finitely often, unless $\eta = 0$. The latter is impossible by the expanding property. □

*Remark.* The convergence of subword frequencies was shown in [7]. We state the associated theorem for the special case of convergence of subletter (i.e. $c \in \mathcal{A}$) frequencies. Due to the previously observed fact that one can use an SMC to define a multitype branching process, the results contained in [2, V.6] can be used to prove the convergence of subletter frequencies (called particles in [2]), while the full statement of [7, Theorem 2.4.10] gives convergence of frequencies for *subwords* (i.e. $w \in \mathcal{A}^n$) as well. Furthermore, [7, Theorem 2.4.10] also provides an algorithm by which one may explicitly calculate such frequencies.

**Theorem 2.5** ([7] Theorem 2.4.10). *For $M$ a matrix indexed by $\mathcal{A} \times \mathcal{A}$ associated to an SMC $(\mathcal{A}^*, (P_c)_{c \in \mathcal{A}})$ with expanding root $a$, defined for $c, c' \in \mathcal{A}$ by*
$$(M)_{c,c'} = \sum_{w \in \mathcal{A}^*} P_{c'}(w) |w|_c$$

if $M$ is primitive (that is, irreducible and aperiodic), then the right eigenvector $e$ of $M$ associated to the eigenvalue 1 with $||e||_1 = 1$ gives the expected frequency of appearance of each $c \in \mathcal{A}$ in the SMC. That is, for $X_n$ the $n^{\text{th}}$ coordinate process for the Markov chain $(\mathcal{A}^*, (P_c)_{c \in \mathcal{A}})$ with initial distribution unit mass on $a$, then for $c \in \mathcal{A}$ and $\mathbb{E}$ the expectation,

$$(e)_c = \lim_{n \to \infty} \mathbb{E} \frac{|X_n|_c}{|X_n|}.$$

We apply this theorem to the second example given below.

**3. Examples.** We gather here a few example SMC's that are particularly analytically tractable.

**3.1. Example 1.** Let $\mathcal{A} = \{a, b\}$, then define $P_a$ and $P_b$ as

$$P_a(w) = \begin{cases} 1/2 & \text{if } w = ab \text{ or } ba \\ 0 & \text{o.w.} \end{cases}$$

$$P_b(w) = \begin{cases} 1 & \text{if } w = b \\ 0 & \text{o.w.} \end{cases}$$

A portion of the SMC started at the symbol $a$ is included in figure 1.

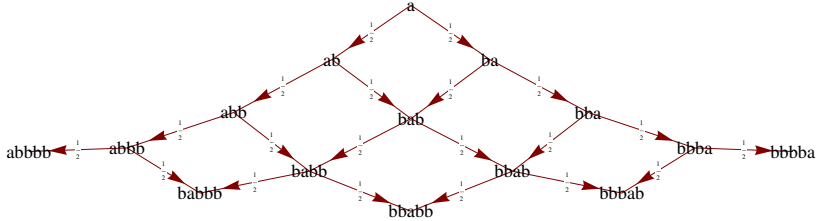

FIGURE 1. A portion of the substitution Markov chain associated to $\Sigma_{\text{eg 1}}$.

A compact way to represent this SMC is to use the following notation:

$$\Sigma_{\text{eg 1}} : \begin{cases} a \to \begin{cases} ab & \text{with prob. } 1/2 \\ ba & \text{with prob. } 1/2 \end{cases} \\ b \to b \end{cases}$$

Note that $a$ is persistent and is a root. This SMC was studied in [7] wherein it was shown that the associated Martin boundary (see section 4) is homeomorphic to the unit interval.

**3.2. Example 2.** For a second example, let $\mathcal{A} = \{a, b\}$, then define $P_a$ and $P_b$ as

$$P_a(w) = \begin{cases} 1/2 & \text{if } w = aa \text{ or } ab \\ 0 & \text{o.w.} \end{cases}$$

$$P_b(w) = \begin{cases} 1 & \text{if } w = ba \\ 0 & \text{o.w.} \end{cases}$$

Written in the alternative notation:

$$\Sigma_{\text{eg 2}} : \begin{cases} a \to \begin{cases} aa \text{ with prob. } 1/2 \\ ab \text{ with prob. } 1/2 \end{cases} \\ b \to ba \end{cases}$$

Note $a$ and $b$ are persistent and the SMC has constant length. We calculate the frequency of appearance of $a, b$ in this SMC via Theorem 2.5. Calculating,

$$M = \begin{array}{c} \\ a \\ b \end{array} \begin{array}{c} a \quad b \\ \begin{bmatrix} 3/2 & 1 \\ 1/2 & 1 \end{bmatrix} \end{array}$$

gives the eigenvector

$$e = \begin{bmatrix} 2/3 \\ 1/3 \end{bmatrix}.$$

Hence, $a$ and $b$ appear with frequency $2/3$ and $1/3$ respectively. This SMC was also studied in [7] wherein it was shown that the associated Martin boundary (see section 4) is homeomorphic to a Cantor space.

Here we add to this by considering a more subtle SMC showing the diversity of this class of Markov chains.

**3.3. Example 3.** Let $\mathcal{A} = \{a, b, c\}$, then define

$$P_a(w) = \begin{cases} 1/4 & \text{if } w = ab, \ ba, \ ac, \text{ or } ca \\ 0 & \text{o.w.} \end{cases}$$

$$P_b(w) = \begin{cases} 1 & \text{if } w = b \\ 0 & \text{o.w.} \end{cases}$$

$$P_c(w) = \begin{cases} 1 & \text{if } w = c \\ 0 & \text{o.w.} \end{cases}$$

Written in the alternative notation:

$$\Sigma_{\text{eg 3}} : \begin{cases} a \to \begin{cases} ab \text{ with prob. } 1/4 \\ ba \text{ with prob. } 1/4 \\ ac \text{ with prob. } 1/4 \\ ca \text{ with prob. } 1/4 \end{cases} \\ b \to b \\ c \to c \end{cases}$$

Note that this SMC is expanding, hence a transient Markov chain. The letter $a$ is persistent and also a root of the Markov chain. This will simplify the calculations in Section 4.

**3.4. Example 4.** Let $\mathcal{A} = \{a, b, c\}$, then define

$$P_a(w) = \begin{cases} q/4 & \text{if } w = ab,\ ba,\ ac,\ \text{or } ca \\ 1-q & \text{if } w = a \\ 0 & \text{o.w.} \end{cases}$$

$$P_b(w) = \begin{cases} 1 & \text{if } w = b \\ 0 & \text{o.w.} \end{cases}$$

$$P_c(w) = \begin{cases} 1 & \text{if } w = c \\ 0 & \text{o.w.} \end{cases}$$

Written in the alternative notation:

$$\Sigma_{eg\ 4}: \begin{cases} a \to \begin{cases} ab \text{ with prob. } q/4 \\ ba \text{ with prob. } q/4 \\ ac \text{ with prob. } q/4 \\ ca \text{ with prob. } q/4 \\ a \text{ with prob. } 1-q \end{cases} \\ b \to b \\ c \to c \end{cases}$$

**3.5. Example 5.** Let $\mathcal{A} = \{a, b, c, d, e\}$, then define

$$P_a(w) = \begin{cases} 1/4 & \text{if } w = ab,\ ba,\ ac,\ \text{or } ca \\ 0 & \text{o.w.} \end{cases}$$

$$P_b(w) = \begin{cases} 1/2 & \text{if } w = bde \text{ or } edb \\ 0 & \text{o.w.} \end{cases}$$

$$P_x(w) = \begin{cases} 1 & \text{if } w = x \text{ and } x \in \{c, d, e\} \\ 0 & \text{o.w.} \end{cases}$$

Written in the alternative notation:

$$\Sigma_{eg\ 5}: \begin{cases} a \to \begin{cases} ab \text{ with prob. } 1/4 \\ ba \text{ with prob. } 1/4 \\ ac \text{ with prob. } 1/4 \\ ca \text{ with prob. } 1/4 \end{cases} \\ b \to \begin{cases} bde \text{ with prob. } 1/2 \\ edb \text{ with prob. } 1/2 \end{cases} \\ c \to c \\ d \to d \\ e \to e \end{cases}$$

The fourth and fifth examples are variations of the third and have Martin kernels which are explicitly computable. However, we focus on the third example for the remainder of this manuscript and, in particular, compute its Martin boundary.

**4. Martin Boundaries.** The Martin boundary is an important topological boundary from a probabilistic and potential theoretic viewpoint as it describes all positive harmonic functions by integrals over this boundary and solves the associated Dirichlet problem. The first thorough treatment of the Martin boundary for Markov chains was by Dynkin [3] who wrote a well-considered and polished account of Hunt's paper [5]. In the context of Markov chains, the standard sources are the books by Kemeny, Snell and Knapp [6], Revuz [13], and Woess [15], [16].

The key difference between the classical literature just mentioned and substitution Markov chains considered here is that a non-degenerate SMC is never irreducible on *any* non-trivial subset of the state-space. Besides Dynkin, the classical literature ([6], [13], [15], [16]) all consider Markov chains $(X, P)$ that are irreducible. Dynkin [3] considered Markov chains $(X, P)$ and initial distributions $\gamma$ referred to by him as *standard measures*: a measure $\gamma : X \to \mathbb{R}$ is called *standard* if for all $y \in X$, $\sum_{x \in X} \gamma(x) \sum_{n=0}^{\infty} P^n(x, y) > 0$. We wish to consider initial distributions that are point masses: $x \in X$, $\gamma = \delta_x$, which fit into that concept if one changes the state space $X$ to be all those states reachable when starting at $a$.

**4.1. Notation.** Here we consider an expanding SMC with root $a$. We first fix the root to be $a \in \mathcal{A}$ and define the *language* to be $\mathcal{L} = \{x \in \mathcal{A}^* : \exists n \geq 0 : P^{(n)}(a, x) > 0 \}$. For $x, y \in \mathcal{L}$, let

$$(4.1) \qquad G(x, y) = \sum_{n=0}^{\infty} P^{(n)}(x, y), \quad K(x, y) = \frac{G(x, y)}{G(a, y)}$$

be the *Green's function* and *Martin kernel* respectively. Choose weights $w_z > 0$ such that $\sum_{z \in \mathcal{L}} \frac{w_z}{G(a, z)} < \infty$. The *Martin metric* defined on $\mathcal{L} \times \mathcal{L}$ is given by

$$(4.2) \qquad \theta(x, y) = \sum_{z \in \mathcal{L}} w_z \left| K(z, x) - K(z, y) \right|.$$

Clearly, $\theta$ is a pseudometric, and if there exist $x \neq y$ such that $\theta(x, y) = 0$, then it must be that for all $z$, $K(z, x) = K(z, y)$. In particular, $K(x, y) = K(x, x) = \frac{1}{G(a,x)} > 0$ and $K(y, x) = K(y, y) = \frac{1}{G(a,y)} > 0$, and hence both $G(x, y) > 0$ and $G(y, x) > 0$, contradicting the expanding property of the SMC.

The *Martin compactification* is given by the completion $\hat{\mathcal{L}}$ of the metric space $(\mathcal{L}, \theta)$, and the *Martin boundary* is given by the compactification minus its interior: $\mathcal{M} = \hat{\mathcal{L}} \setminus \mathcal{L}$. Equivalently (see [7]), for $x_n, y_n \in \mathcal{L}$, the Martin boundary is given by equivalence classes of sequences $(x_n)_{n \geq 0}$, $(y_n)_{n \geq 0}$ which leave every finite subset, such that $(x_n)_{n \geq 0}$ and $(y_n)_{n \geq 0}$ are equivalent if and only if for all $z \in \mathcal{L}$, $\lim_{n \to \infty} K(z, x_n) = \lim_{n \to \infty} K(z, y_n)$.

**Lemma 4.1.** *A constant length, persistent substitution Markov chain $(\mathcal{A}^*, (P_c)_{c \in \mathcal{A}})$ with an expanding root $a \in \mathcal{A}$ that satisfies $|\mathrm{supp}(P_a)| \geq 2$ and $\forall c \in \mathcal{A}$, $\mathrm{supp}(P_c) \setminus \bigcup_{c \neq b \in \mathcal{A}} \mathrm{supp}(P_b) \neq \emptyset$ has a non-trivial Martin boundary.*

*Proof.* We will demonstrate this fact by explicitly constructing a strictly positive, non-constant harmonic function. Given such an SMC $(\mathcal{A}^*, (P_c)_{c \in \mathcal{A}})$, define a deterministic substitution $\sigma$: for each $c \in \mathcal{A}$, fix a $w_c \in \text{supp}(P_c) \setminus \bigcup_{c \neq b \in \mathcal{A}} \text{supp}(P_b)$ and let $\sigma(c) = w_c$. Now since the SMC has an expanding root $a$, we have $|\sigma^n(a)|_a = 1$ and $|\sigma^{n+1}(a)| > |\sigma^n(a)|$. Now let $(s_n)_{n \geq 0}$ be a sequence of real numbers defined by the recursion relation: $s_0 = 1$, $s_{n+1} = \frac{s_n - k}{P(\sigma^n(a), \sigma^{n+1}(a))} + k$ for a fixed $k \in \mathbb{R}^+$, $k \neq 1$. Then for $w \in \mathcal{A}^*$, define

$$f(w) = \begin{cases} s_n, & \text{if } w = \sigma^n(a) \\ k, & \text{otherwise} \end{cases}.$$

Then for $w \in \mathcal{A}^*$ with $w \neq \sigma^n(a)$ for any $n$, since the SMC is constant length, if $P(w, v) > 0$, then $v \neq \sigma^m(a)$ for any $m$. Hence we have

$$\begin{aligned} Pf(w) &= \sum_{v \in \mathcal{A}^*} P(w, v) f(v) \\ &= \sum_{v \in \mathcal{A}^*} P(w, v) k \\ &= k \\ &= f(w). \end{aligned}$$

Now in the case $w = \sigma^n(a)$ for some $n$, we have

$$\begin{aligned} Pf(w) &= \sum_{v \neq \sigma^{n+1}(a)} P(\sigma^n(a), v) f(v) + P(\sigma^n(a), \sigma^{n+1}(a)) f(\sigma^{n+1}(a)) \\ &= k(1 - P(\sigma^n(a), \sigma^{n+1}(a))) + P(\sigma^n(a), \sigma^{n+1}(a)) f(\sigma^{n+1}(a)) \\ &= k + P(\sigma^n(a), \sigma^{n+1}(a))(f(\sigma^{n+1}(a)) - k) \\ &= k + P(\sigma^n(a), \sigma^{n+1}(a))(s_{n+1} - k) \\ &= f(w). \end{aligned}$$

So in either case, $Pf(w) = f(w)$ so $f$ is harmonic. □

**4.2. Martin Boundary of Example 3.** Let $\mathcal{L} \subseteq \{a, b, c\}^*$ be the language (or state space) of the substitution Markov chain $\Sigma_{eg\ 3}$. Note that for each $w \in \mathcal{L}$, $|w|_a = 1$, and each application of the random substitution adds a $b$ or $c$ immediately to the left or right of the $a$. We will take advantage of this fact by treating the symbol $a$ as an "origin". For $w, v \in \mathcal{A}^*$, let

(4.3) $$d_s(w, v) = 2^{-\max\{n:\ w_1 \ldots w_n = v_1 \ldots v_n\}}$$

be the standard word metric. We also denote the reverse of the word $w = w_1 \ldots w_n$ by $\bar{w} = w_n \ldots w_1$.

First of all notice that we can calculate the Martin kernel explicitly: Let $z = HaT$ be a fixed word in $\mathcal{L}$ where the lengths of the subwords $H$ and $T$ are $h$ and $t$. We think of $H$ and $T$ as being the "head" and "tail" respectively of $z$. Then one can easily prove

**Lemma 4.2.**

$$G(z, x) = \begin{cases} 4^{-l_x - r_x} \binom{l_x + r_x}{r_x} & \text{if } x = HL_x a R_x T, |L_x| = l_x; |R_x| = r_x \\ 0 & \text{otherwise} \end{cases}$$

and
$$K(z,x) = \begin{cases} 4^{h+t} \frac{(h+l_x)_h (t+r_x)_t}{(h+t+l_x+r_x)_{h+t}} & \text{if } x = HL_x a R_x T; |L_x| = l_x; |R_x| = r_x \\ 0 & \text{otherwise.} \end{cases}$$

Here we use the notation
$$(n)_k = n(n-1)...(n-k+1).$$

*Proof.* By definition, $G(HaT, x) > 0$ if and only if $x = HL_x a R_x T$ for some finite words $L_x$ and $R_x$ of length $l_x$ and $r_x$, and

$$G(HaT, x) = P^{(n)}(HaT, x) = \begin{cases} 4^{-l_x-r_x} \binom{l_x+r_x}{r_x} & \text{if } n = l_x + r_x \\ 0 & \text{otherwise.} \end{cases}$$

Then, omitting the index $x$,

$$K(HaT, HLaRT) = \frac{G(HaT, HLaRT)}{G(a, HLaRT)} = \frac{4^{-l-r}\binom{l+r}{r}}{4^{-l-r-h-t}\binom{l+r+h+t}{h+l}}$$
$$= 4^{h+t} \frac{(l+r)!(h+l)\cdots(l+1)l!(t+r)\cdots(r+1)r!}{l!r!(l+r+h+t)\cdots(l+r+1)(l+r)!}.$$
$\square$

Let $(x_n)_{n\geq 1} \in \mathcal{L}^\mathbb{N}$ be a Cauchy sequence in the metric $\theta$ from definition 4.2, called a $\theta$-Cauchy sequence in the sequel. Then for $z = HaT \in \mathcal{L}$, $(K(z, x_n))_{n\geq 0}$ is Cauchy in $\mathbb{R}$, whence
$$\lim_{n\to\infty} K(z, x_n)$$
exists and this limit equals zero or
$$4^{h+t} \lim_{n\to\infty} \left(\frac{l_{x_n}}{l_{x_n}+r_{x_n}}\right)^h \left(\frac{r_{x_n}}{l_{x_n}+r_{x_n}}\right)^t,$$
in case that $x_n = HL_{x_n} a R_{x_n} T$ for all $n$ large enough, where $l_{x_n} = |L_{x_n}|$ and $r_{x_n} = |R_{x_n}|$.

Given two $\theta$-Cauchy sequences $(x_n)_{n\geq 0}$ and $(y_n)_{n\geq 0}$ we write
$$x_n = L_{x_n} a R_{x_n}$$
$$y_n = L_{y_n} a R_{y_n}$$
with $L_{x_n}, L_{y_n} \in \{b, c\}^* \cup \{\epsilon\}$ (where $\epsilon$ is the empty word) denoting the portions of $x_n$ and $y_n$ to the left of the "origin" $a$, and similarly for $R_{x_n}, R_{y_n}$ to the right.

Now we can define an equivalence relation " $\sim$ " on the set of Cauchy sequences on the state space $\mathcal{L}$. Let $(x_n)_{n\geq 0} \sim (y_n)_{n\geq 0}$ if the following conditions are satisfied:

**I:** $\lim_{n\to\infty} \frac{|L_{x_n}|}{|L_{x_n}|+|R_{x_n}|} = \lim_{n\to\infty} \frac{|L_{y_n}|}{|L_{y_n}|+|R_{y_n}|} =: \lambda$.
**II(a):** If $\lambda \in (0,1)$ then $\lim_{n\to\infty} d_s(L_{x_n}, L_{y_n}) = \lim_{n\to\infty} d_s(\overleftarrow{R}_{x_n}, \overleftarrow{R}_{y_n}) = 0$.
**II(b):** If $\lambda = 0$ then $\lim_{n\to\infty} d_s(\overleftarrow{R}_{x_n}, \overleftarrow{R}_{y_n}) = 0$.
**II(c):** If $\lambda = 1$ then $\lim_{n\to\infty} d_s(L_{x_n}, L_{y_n}) = 0$.

We aim to show that $\mathcal{L}^{\mathbb{N}}/\sim$ naturally gives rise to the Martin boundary $\mathcal{M}$.

We are now in a position to demonstrate

**Lemma 4.3.** *The following two conditions are equivalent for two $\theta$-Cauchy sequences $(x_n)_{n\geq 0}$ and $(y_n)_{n\geq 0}$:*
*(1) $(x_n)_{n\geq 0} \sim (y_n)_{n\geq 0}$*
*(2) $(x_n)_{n\geq 0}$ and $(y_n)_{n\geq 0}$ are equivalent $\theta$-Cauchy sequences.*

*Proof.* (1)$\Longrightarrow$(2). Let $(x_n)_{n\geq 0} \sim (y_n)_{n\geq 0}$ be $\theta$-Cauchy sequences. Write $x_n = L_{x_n}aR_{x_n}$ and $y_n = L_{y_n}aR_{y_n}$ with $R_{x_n} = (R_{x_n}(1),...,R_{x_n}(r_{x_n}))$, $R_{y_n} = (R_{y_n}(1),...,R_{y_n}(r_{y_n}))$, $L_{x_n} = (L_{x_n}(1),...,L_{x_n}(l_{x_n}))$ and $L_{y_n} = (L_{y_n}(1),...,L_{y_n}(l_{y_n}))$. Let $z = HaT \in \mathcal{L}$ with $h = |H|$ and $t = |T|$. By **II** there exists $n_0 \geq 0$ such that for $n \geq n_0$ we have the following: in the case of **II(a)**, that

$$(L_{x_n}(1),...,L_{x_n}(h)) = (L_{y_n}(1),...,L_{y_n}(h)) \quad \text{and}$$
$$(R_{x_n}(r_{x_n}),...,R_{x_n}(r_{x_n}-t+1)) = (R_{y_n}(r_{y_n}),...,R_{y_n}(r_{y_n}-t+1)).$$

In the case of **II(b)**, that

$$(R_{x_n}(r_{x_n}),...,R_{x_n}(r_{x_n}-t+1)) = (R_{y_n}(r_{y_n}),...,R_{y_n}(r_{y_n}-t+1)).$$

And in the case of **II(c)**, that

$$(L_{x_n}(1),...,L_{x_n}(h)) = (L_{y_n}(1),...,L_{y_n}(h)).$$

In case **II(a)** it follows that either $G(z,x_n) = G(z,y_n) = 0$ for all $n \geq n_0$, or else by Lemma 4.2

$$\lim_{n\to\infty} K(z,x_n) = \lim_{n\to\infty} 4^{h+t} \left(\frac{|L_{x_n}|}{|L_{x_n}|+|R_{x_n}|}\right)^h \left(\frac{|R_{x_n}|}{|L_{x_n}|+|R_{x_n}|}\right)^t + o(1)$$
$$= 4^{h+t}\lambda^h(1-\lambda)^t = \lim_{n\to\infty} K(z,y_n).$$

In case **II(b)** it follows for $h > 0$ by the same lemma that

$$\lim_{n\to\infty} K(z,x_n) = 0 = \lim_{n\to\infty} K(z,y_n).$$

For $h = 0$ we either have $G(z,x_n) = G(z,y_n) = 0$ for all $n \geq n_0$, or else by Lemma 4.2

$$\lim_{n\to\infty} K(z,x_n) = \lim_{n\to\infty} 4^t \left(\frac{|R_{x_n}|}{|L_{x_n}|+|R_{x_n}|}\right)^t + o(1)$$
$$= 4^t(1-\lambda)^t = \lim_{n\to\infty} K(z,y_n).$$

The case **II(c)** is shown in the same way.

(2)$\Longrightarrow$(1). This is easily shown using similar arguments as in the first part of the proof. $\square$

We have thus shown that the Martin boundary is characterized by $\sim$. We use this representation of equivalence classes of $\theta$-Cauchy sequences to describe the Martin boundary as a metric space.

**Proposition 4.4.** *The Martin boundary $\mathcal{M}$ equipped with the metric*

$$\theta(\xi, \eta) = \sum_{z \in \mathcal{L}} 16 \cdot 4^{-2|z|} |K(z, \xi) - K(z, \eta)|,$$

*where $K(z, \xi) = \lim_{n \to \infty} K(z, x_n)$ with $\xi$ being represented by any of its $\theta$-Cauchy sequences $(x_n)_{n \geq 0}$, is isometric to the space*

$$B = \left\{ (\lambda, \mathbf{x}) : \lambda \in [0, 1]; \mathbf{x} = \begin{cases} (L_\mathbf{x}, R_\mathbf{x}) \in \{b, c\}^\mathbb{N} \times \{b, c\}^\mathbb{N} & \text{if } 0 < \lambda < 1 \\ R_\mathbf{x} \in \{b, c\}^\mathbb{N} & \text{if } \lambda = 0 \\ L_\mathbf{x} \in \{b, c\}^\mathbb{N} & \text{if } \lambda = 1 \end{cases} \right\}$$

*with the metric*

$$\rho((\lambda, \mathbf{x}), (\mu, \mathbf{y})) = \sum_{h=0}^{n} \sum_{t=0}^{m} 4^{-t-h} |\lambda^h (1-\lambda)^t - \mu^h (1-\mu)^t|$$
$$+ \frac{16}{(3+\lambda)(4-\lambda)} \left[ \left(\frac{\lambda}{4}\right)^{n+1} + \left(\frac{1-\lambda}{4}\right)^{m+1} - \left(\frac{\lambda}{4}\right)^{n+1} \left(\frac{1-\lambda}{4}\right)^{m+1} \right]$$
$$+ \frac{16}{(3+\mu)(4-\mu)} \left[ \left(\frac{\mu}{4}\right)^{n+1} + \left(\frac{1-\mu}{4}\right)^{m+1} - \left(\frac{\mu}{4}\right)^{n+1} \left(\frac{1-\mu}{4}\right)^{m+1} \right]$$

*where*
$$n = \max\{k : R_\mathbf{x}(k) = R_\mathbf{y}(k)\} \qquad m = \min\{k : L_\mathbf{x}(k) = L_\mathbf{y}(k)\}$$

*Proof.* Let $\xi$ be an equivalence class of $\theta$-Cauchy sequences. By Lemma 4.3 such a sequence is uniquely characterized by a map $F : \mathcal{M} \to B$ defined by any representative $x_n = L_n a R_n$ from $\xi$ via
$$F(\xi) = (\lambda, \mathbf{x})$$

with $\lambda = \lim_{n \to \infty} \frac{|L_n|}{|L_n| + |R_n|}$ and $\mathbf{x} = (\lim_{n \to \infty} L_n, \lim_{n \to \infty} R_n)$ in case of $0 < \lambda < 1$, $\mathbf{x} = \lim_{n \to \infty} R_n$ in case of $\lambda = 0$, and $\mathbf{x} = \lim_{n \to \infty} L_n$ in case of $\lambda = 1$.

Thus it is sufficient to calculate the distance between $\xi$ and $\eta$ using their representations and to express it in terms of $F(\xi)$ and $F(\eta)$.

Let $\xi$ and $\eta$ be two equivalence classes of $\theta$-Cauchy sequences represented by $F(\xi) = (\lambda, \mathbf{x})$ and $F(\eta) = (\mu, \mathbf{y})$. Assume first that $\lambda, \mu \in (0, 1)$. Write $\mathbf{x} = (L_x(1), L_x(2), ...; R_x(1), R_x(2), ...)$ and $\mathbf{y} = (L_y(1), L_y(2), ...; R_y(1), R_y(2), ...)$. Let $n$ and $m$ be as in the statement of the proposition. For each $h, t \geq 0$ there is exactly one $z_x(h, t) = z = HaT$ with $|H| = h$, $|T| = t$ and $0 < K(z, \xi) = 4^{h+t} \lambda^h (1-\lambda)^t$ and likewise for $\eta$. For all $h \leq n$,

$t \leq m$, $z_x(h,t) = z_y(h,t)$. Hence, counting $z$'s in terms of $h$ and $t$,

$$\begin{aligned}
\theta(\xi,\eta) &= \sum_{h=0}^{\infty}\sum_{t=0}^{\infty} 4^{-2t-2h}\{|K(z_x(h,t),\xi) - K(z_x(h,t),\eta)| \\
&\quad + \mathbb{I}_{z_x(h,t)\neq z_y(h,t)}|K(z_y(h,t),\xi) - K(z_y(h,t),\eta)|\} \\
&= \sum_{h=0}^{n}\sum_{t=0}^{m} 4^{-t-h}|\lambda^h(1-\lambda)^t - \mu^h(1-\mu)^t| \\
&\quad + \sum_{h=n+1}^{\infty}\sum_{t=0}^{\infty} 4^{-t-h}\lambda^h(1-\lambda)^t + \sum_{h=0}^{n}\sum_{t=m+1}^{\infty} 4^{-t-h}\lambda^h(1-\lambda)^t \\
&\quad + \sum_{h=n+1}^{\infty}\sum_{t=0}^{\infty} 4^{-t-h}\mu^h(1-\mu)^t + \sum_{h=0}^{n}\sum_{t=m+1}^{\infty} 4^{-t-h}\mu^h(1-\mu)^t \\
&= \sum_{h=0}^{n}\sum_{t=0}^{m} 4^{-t-h}|\lambda^h(1-\lambda)^t - \mu^h(1-\mu)^t| \\
&\quad + \frac{16}{(3+\lambda)(4-\lambda)}\left[\left(\frac{\lambda}{4}\right)^{n+1} + \left(\frac{1-\lambda}{4}\right)^{m+1} - \left(\frac{\lambda}{4}\right)^{n+1}\left(\frac{1-\lambda}{4}\right)^{m+1}\right] \\
&\quad + \frac{16}{(3+\mu)(4-\mu)}\left[\left(\frac{\mu}{4}\right)^{n+1} + \left(\frac{1-\mu}{4}\right)^{m+1} - \left(\frac{\mu}{4}\right)^{n+1}\left(\frac{1-\mu}{4}\right)^{m+1}\right]
\end{aligned}$$

$\square$

Recall that the Cantor space is defined by $C = \{0,1\}^{\mathbb{N}} = \{(u^i)_{i\geq 1} : u^i \in \{0,1\}\}$ equipped with the metric $d_s^r$ for any $r > 0$ where $d_s$ is defined in (4.3).

We are now in the position to demonstrate the following theorem.

**Theorem 4.5.** *The Martin boundary of the SMC $\Sigma_{\text{eg 3}}$ equipped with the Martin metric is homeomorphic to the product space of the unit interval with the Cantor space, equipped with the metric*

$$d((\lambda,u),(\mu,v)) = |\lambda - \mu| + d_s(u,v).$$

*The homeomorphism can be chosen so that its inverse is a Lipschitz continuous map. Moreover both spaces are not Lipschitz equivalent for any of the metrics $|\lambda - \mu| + d_s^r(u,v)$ where $r > 0$.*

*Proof.* We first define a map

$$\Phi : B \to [0,1] \times C.$$

To this end let $\lambda \in [0, \frac{1}{2}]$ and define

$$p_{k-1} = \lfloor -k\log_2(\lambda/4)\rfloor \qquad k \in \{1,2,3,...\}.$$

For $\lambda = 0$ this means that no sequence $p_k$ is chosen. Then define for $\lambda \neq 0$

$$\Phi(\lambda,(l_1,l_2,...;r_1,r_2,...)) = (\lambda,(u^k)_{k\geq 0})$$

$u^{p_k} = \mathbb{I}_{r_k=c}$ and $(u^i)_{i \notin \{p_k : k \geq 0\}} = (\mathbb{I}_{l_1=c}, \mathbb{I}_{l_2=c}, ...)$. For $\lambda = 0$

$$\Phi(0, (r_1, r_2, ...)) = (0, (\mathbb{I}_{r_1=c}, \mathbb{I}_{r_2=c}, ...)).$$

In case $\lambda \in (\frac{1}{2}, 1]$ define $\Phi$ in the same way replacing $\lambda$ by $1 - \lambda$, and switching $r_i$ with $l_i$.

We want to show that $\Phi^{-1}$ is Lipschitz continuous with respect to the metric $\rho$ on $B$ and $d$ on $[0,1] \times C$. Then $\Phi$ is also continuous, since the domain and the image of $\Phi$ are compact.

To show this part of the theorem we show that there is a constant $\kappa > 0$ such that for all $(\lambda, \mathbf{x}), (\mu, \mathbf{y}) \in B$

$$\rho((\lambda, \mathbf{x}), (\mu, \mathbf{y})) \leq \kappa d(\Phi((\lambda, \mathbf{x})), \Phi((\mu, \mathbf{y}))).$$

We consider the case $\lambda, \mu \leq \frac{1}{2}$ only since the other cases are similar. Use $m$ and $n$ as in proposition 4.4. Now, note that the derivative of $x^h(1-x)^t$ is bounded by $h+t$ on the unit interval. Hence by the mean value theorem and finiteness of the expectation of a geometric distribution, we find a constant $\kappa_1$ such that

$$\sum_{h=0}^{n} \sum_{t=0}^{m} 4^{-(h+t)} |\lambda^h(1-\lambda)^t - \mu^h(1-\mu)^t| \leq \kappa_1 |\lambda - \mu|.$$

By definition of the map $\Phi$, for $\Phi((\lambda, \mathbf{x})) = (\lambda, (u^i)_{i\geq 0})$ and $\Phi((\mu, \mathbf{y})) = (\mu, (v^i)_{i\geq 0})$ we have that $u^i \neq v^i$ for the first index $i$ if and only if $i = p_{n+1}$ or for $i = m + k + 1$ (or for $i = m + k + 2$) where $p_k \leq m < p_{k+1}$ for some $k \leq n$. Then

$$\left(\frac{\lambda}{4}\right)^{n+1} = 2^{(n+1)\log_2\left(\frac{\lambda}{4}\right)} \leq 2 \cdot 2^{-p_{n+1}}$$

and, since $k \leq m$,

$$\begin{aligned}\left(\frac{1-\lambda}{4}\right)^{m+1} &= 2^{(m+1)\log_2\left(\frac{1-\lambda}{4}\right)} \\ &= 2 \cdot 2^{-(m+k+2)} \cdot 2^{(m+1)\left(1+\log_2\left(\frac{1-\lambda}{4}\right)\right)+k} \\ &\leq 2 \cdot 2^{-(m+k+2)} \cdot 2^{(m+1)\left(2+\log_2\left(\frac{1-\lambda}{4}\right)\right)}\end{aligned}$$

Since $2 + \log_2\left(\frac{1-\lambda}{4}\right) \leq 0$, the last expression is $O\left(2^{-(m+k+1)}\right)$. Hence there exists a constant $\kappa_2$ such that

$$\frac{16}{(3+\lambda)(4-\lambda)}\left[\left(\frac{\lambda}{4}\right)^{n+1} + \left(\frac{1-\lambda}{4}\right)^{m+1}\right] \leq \kappa_2 d_s((u^i)_{i\geq 0}, (v^i)_{i\geq 0}).$$

Putting this and similar estimates for the other cases together shows the claim.

It is left to show that there cannot be a bi-Lipschitz map between both metric spaces. First note that $[0,1] \times C$ equipped with any metric $|\lambda - \mu| + d_s^r(u, v)$ has the property that for any non-empty open set the Hausdorff dimension of this set equals the Hausdorff dimension of the space $[0,1] \times C$. The Hausdorff dimension of the sets

$$U(\lambda, \eta) = \{(\mu, \mathbf{x}) \in B : \lambda - \eta < \mu < \lambda + \eta\}$$

is not constant and varies with $\lambda$ for small $\eta > 0$. If $\Phi$ is a bi-Lipschitz continuous homeomorphism its image $\Phi(U(\lambda, \eta))$ is open and has the same Hausdorff dimension as $U(\lambda, \eta)$, which is a contradiction. □

**Theorem 4.6.** *The Martin boundary of example $\Sigma_{\text{eg 3}}$ is bi-Lipschitz equivalent to the space*

$$D = \left\{ \left(\lambda, \sum_{k=1}^{\infty} x_k \left(\frac{\lambda}{4}\right)^k, \sum_{k=1}^{\infty} y_k \left(\frac{1-\lambda}{4}\right)^k \right) : 0 \leq \lambda \leq 1,\ x_k, y_k \in \{1,3\}\ \text{for } k \geq 1 \right\}$$

*equipped with the Euclidean metric.*

*Proof.* We use the Euclidean metric $\|\cdot - \cdot\|_2$ in $\mathbb{R}^3$ and the map

$$\Psi : B \to D$$

$$(\lambda, (x_k), (y_k)) \mapsto \left(\lambda, \sum_{k=1}^{\infty} (1 + 2\mathbb{I}_{x_k=c}) \left(\frac{\lambda}{4}\right)^k, \sum_{k=1}^{\infty} (1 + 2\mathbb{I}_{y_k=c}) \left(\frac{1-\lambda}{4}\right)^k \right).$$

Let $(\lambda, (x_k), (y_k)), (\mu, (u_k), (v_k)) \in B$. To be more precise $(\lambda, \mathbf{x}) \in B$ with $L_\mathbf{x} = x_k$ and $R_\mathbf{y} = y_k$ and similarly with $(\mu, \mathbf{u})$. Let $\tilde{x}_k = 1 + 2\mathbb{I}_{x_k=c}$, $\tilde{y}_k = 1 + 2\mathbb{I}_{y_k=c}$, $\tilde{u}_k = 1 + 2\mathbb{I}_{u_k=c}$, $\tilde{v}_k = 1 + 2\mathbb{I}_{v_k=c}$ and let $m$ and $n$ be as in Proposition 4.4.

Using the mean value theorem, for some constants $C_1$ and $C_2$ we obtain

$$\left|\sum_{k=1}^{\infty} \left(\tilde{x}_k \left(\frac{\lambda}{4}\right)^k - \tilde{u}_k \left(\frac{\mu}{4}\right)^k\right)\right| \leq 3 \sum_{k=n+1}^{\infty} \left(\frac{\lambda}{4}\right)^k + \left(\frac{\mu}{4}\right)^k + \left|\sum_{k=1}^{n} \tilde{x}_k 4^{-k}(\lambda^k - \mu^k)\right|$$

$$\leq C_1 \left[\left(\frac{\lambda}{4}\right)^{n+1} + \left(\frac{\mu}{4}\right)^{n+1}\right] + 3|\lambda - \mu| \sum_{k=1}^{n} 4^{-k} \left|\frac{\lambda^k - \mu^k}{\lambda - \mu}\right|$$

$$\leq C_2 \rho((\lambda, (x_k), (y_k)), (\mu, (u_k), (v_k))).$$

It follows from this and the similar estimate replacing $\lambda$ by $1 - \lambda$ and $\mu$ by $1 - \mu$ that there is a constant $C$ such that

$$\left\|\left(\lambda, \sum_{k=1}^{\infty} \tilde{x}_k \left(\frac{\lambda}{4}\right)^k, \sum_{k=1}^{\infty} \tilde{y}_k \left(\frac{1-\lambda}{4}\right)^k\right) - \left(\mu, \sum_{k=1}^{\infty} \tilde{u}_k \left(\frac{\mu}{4}\right)^k, \sum_{k=1}^{\infty} \tilde{v}_k \left(\frac{1-\mu}{4}\right)^k\right)\right\|_2$$

$$\leq C\rho((\lambda, (x_k), (y_k)), (\mu, (u_k), (v_k))).$$

A lower bound is obtained as follows. In the case that

$$|\lambda - \mu| \geq \frac{1}{2} \max\left\{\left(\frac{\lambda}{4}\right)^{n+1}, \left(\frac{\mu}{4}\right)^{n+1}, \left(\frac{1-\lambda}{4}\right)^{m+1}, \left(\frac{1-\mu}{4}\right)^{m+1}\right\},$$

we have

$$\left\|\left(\lambda, \sum_{k=1}^{\infty} \tilde{x}_k \left(\frac{\lambda}{4}\right)^k, \sum_{k=1}^{\infty} \tilde{y}_k \left(\frac{1-\lambda}{4}\right)^k\right) - \left(\mu, \sum_{k=1}^{\infty} \tilde{u}_k \left(\frac{\mu}{4}\right)^k, \sum_{k=1}^{\infty} \tilde{v}_k \left(\frac{1-\mu}{4}\right)^k\right)\right\|_2$$

$$\geq |\lambda - \mu| \geq C_3 \rho((\lambda, (x_k), (y_k)), (\mu, (u_k), (v_k)))$$

for some constant $C_3 > 0$.

In the case of

$$\left(\frac{\lambda}{4}\right)^{n+1} \geq \max\left\{\left(\frac{\mu}{4}\right)^{n+1}, \left(\frac{1-\lambda}{4}\right)^{m+1}, \left(\frac{1-\mu}{4}\right)^{m+1}, 2|\lambda - \mu|\right\}$$

we get for some constant $C_4 > 0$ that

$$\left\|\left(\lambda, \sum_{k=1}^{\infty} \tilde{x}_k \left(\frac{\lambda}{4}\right)^k, \sum_{k=1}^{\infty} \tilde{y}_k \left(\frac{1-\lambda}{4}\right)^k\right) - \left(\mu, \sum_{k=1}^{\infty} \tilde{u}_k \left(\frac{\mu}{4}\right)^k, \sum_{k=1}^{\infty} \tilde{v}_k \left(\frac{1-\mu}{4}\right)^k\right)\right\|_2$$

$$\geq \left|\sum_{k=1}^{\infty} 4^{-k}(\tilde{x}_k \lambda^k - \tilde{u}_k \mu^k)\right|$$

$$= \left|(\tilde{x}_{n+1} - \tilde{u}_{n+1})4^{-n-1}\lambda^{n+1} + \sum_{k=n+2}^{\infty} (\tilde{x}_k - \tilde{u}_k)4^{-k}\lambda^k + \sum_{k=1}^{\infty} \tilde{u}_k 4^{-k}(\lambda^k - \mu^k)\right|$$

$$\geq 2\left(\frac{\lambda}{4}\right)^{n+1} - 8\frac{4^{-n-2}\lambda^{n+2}}{4-\lambda} - 3|\lambda - \mu|\sum_{k=1}^{\infty} k 4^{-k}$$

$$= \left(\frac{\lambda}{4}\right)^{n+1}\left(2 - \frac{2\lambda}{4-\lambda}\right) - \frac{4}{3}|\lambda - \mu|$$

$$\geq \frac{4}{3}\left(\frac{\lambda}{4}\right)^{n+1} - \frac{4}{3}|\lambda - \mu|$$

$$\geq C_4 \rho((\lambda, (x_k), (y_k)), (\mu, (u_k), (v_k))).$$

The other 3 cases are treated similarly. This ends the proof. □

Included in figure 2 is an approximation of the space $D$ defined in theorem 4.6.

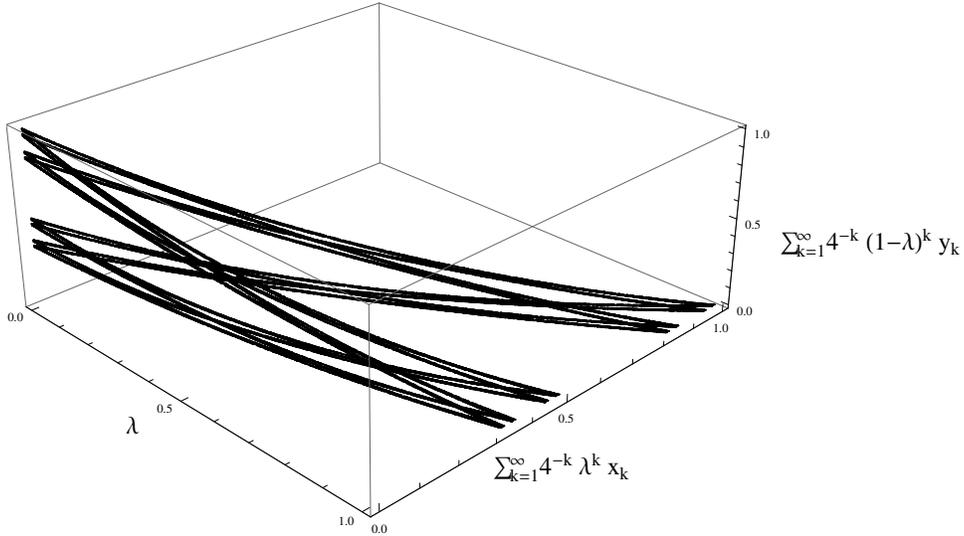

FIGURE 2. An approximate representation of the space $D$.

**References.**

1. G. Alsmeyer and U. Rösler. The martin entrance boundary of the galton-watson process. *Ann. Inst. H. Poincar Probab. Statist.*, 42(5):591–606, 2006.

2. K. B. Athreya and P. E. Ney. *Branching Processes*. Springer-Verlag, Berlin, 1972.
3. E. B. Dynkin. Boundary theory of Markov processes (the discrete case). *Uspekhi Mat. Nauk*, 24(2):1–42, 1969.
4. S. Evans, R. Grübler, and A. Wakolbinger. Trickle-down processes and their boundaries. *Electron. J. Probab.*, 17(1):58pp, 2012.
5. G. Hunt. Markoff chains and Martin boundaries. *Illinois J. Math.*, 4(3):313–340, 1960.
6. J. Kemeny, J. Snell, and A. Knapp. *Denumerable Markov Chains*. Springer-Verlag, New York, NY., 2nd edition, 1965.
7. D. Koslicki. *Substitution Markov chains with applications to molecular evolution*. PhD thesis, Pennsylvania State University, State College, PA, USA, 2012.
8. J. Peyrière. Random beadsets and birth processes with interaction. *IBM Research Report RC-7417*, pages 1–19, 1978.
9. J. Peyrière. Process of birth interaction with the neighbors, evolution of charts. *Ann. Inst. Fourier. Grenoble*, 31(4):187–218, 1981.
10. J. Peyriére. Substitutions aléatoires itérées. *Seminaire de Théorie des Nombres de Bordeaux*, 17:1–9, 1981.
11. J. Peyrière. Frequence des motifs dans les suites doubles invariantes par une subsitution. *Ann. sc. math. Québec*, 11(1):133–138, 1987.
12. M. Queffélec. *Substitution Dynamical Systems - Spectral Analysis, Lecture Notes in Mathematics, 1294*. Springer-Verlag, Berlin, 1980.
13. D. Revuz. *Markov Chains*. North-Holland, Amsterdam, 1975.
14. D. J. Wing. *Notions of Complexity in Substitution Dynamical Systems*. PhD thesis, Oregon State University, Corvallis, OR, USA, 2011.
15. W. Woess. *Random Walks on Infinite Graphs and Groups*. Cambridge University Press, Cambridge, 2000.
16. W. Woess. *Denumerable Markov Chains*. European Math. Soc. Publishing House, Zürich, Switzerland, 2009.

Department of Mathematics, Oregon State University, Corvallis, OR 97331
**Email address**: **david.koslicki@math.oregonstate.edu**

Mathematics Department, The Pennsylvania State University, State College, PA 16802
**Email address**: **denker@math.psu.edu**